\documentclass[12pt]{article}
\usepackage{amsmath}
\usepackage{amssymb}
\usepackage{amscd}
\newtheorem{lemme}{Lemme}
\newtheorem{theoreme}{Th\'eor\`eme}
\newtheorem{proposition}{Proposition}
\newtheorem{corollaire}{Corollaire}

\newenvironment{exemple}
  {\begin{trivlist}\item[]\textbf{Exemple}}{\end{trivlist}}
\newcounter{numeroremarque}
\newenvironment{remarque}
  {\addtocounter{numeroremarque}{1}
\begin{trivlist}\item[]\textbf{Remarque \thenumeroremarque}}{\end{trivlist}}
\newcounter{numerodefinition}

\newcounter{numeroquestion}

\newcounter{numeroconjecture}

\newenvironment{preuve}{\begin{trivlist}\item[]\textit{Preuve.}}
{\item[] \hfill $\square$\end{trivlist}}
\begin{document}
\title{Remarque sur les fonctions ayant le m\^eme ensemble de Julia
\footnote{{\bf Classification math\'ematique:} 30D05, 58F23.
\ \ \ \ \ \ \ \ \ \ \ \ \ \ \ \ \ \ \ \ \ \ \ \ \ \ \
\break
{\bf Mots cl\'es:} ensemble de Julia, 
ensemble analytiquement lamin\'e, suite d'it\'er\'es.}
}
\author{{Tien-Cuong} DINH\footnote{
Math\'ematique-B\^atiment 425,
Universit\'e Paris-Sud, 91405 ORSAY Cedex (France).
}}
\maketitle
%
%
\begin{abstract} \begin{center}
{\bf On the rationnal functions having the same Julia set}
\end{center}
Let $f$ and $g$ two rationnal functions having the same Julia set
$J_f$. Lets suppose that $f$ has a rational indifferent periodic point
and that the critical set of $f$ is disjoint of $J_f$. Then or $J_f$
has to be equal to $\mathbb{P}^1$, a circle, an arc of a
circle for some coordinate or $f$ and $g$ has to verify an equation of
the type:
$$f^{m_1}\circ g\circ f^{m_2}\circ g\circ\cdots\circ f^{m_n}\circ g=f^m.$$ 
\end{abstract}
\section{Introduction}
On note $\mathbb{P}^1:=\mathbb{C}\cup\{\infty\}$
  la droite projective complexe. On appelle  
  {\it fonction  rationnelle} tout quotient $f=P/Q$ o\`u $P$ et
  $Q$ sont deux polyn\^omes sans facteur commun. {\it Le degr\'e} de $f$
  est le maximum des degr\'es de $P$ et de $Q$. La fonction $f$ d\'efinit un
  endomorphisme holomorphe de $\mathbb{P}^1$. Par la
  suite, on note $fg:=f\circ g$, $f^k:=f\circ\cdots\circ f$ ($k$
  fois) et $f^0:=\mbox{Id}$ 
  pour toutes fonctions rationnelles $f$ et $g$. Supposons que $\deg
  f\geq 2$ et $\deg g\geq 2$.\par
On appelle
  {\it l'ensemble de Fatou} de $f$ l'ouvert maximal $F_f$ o\`u la
  suite $\{f^k\}_{k=1}^\infty$ est \'equicontinue. {\it L'ensemble de
  Julia} $J_f$ est le compl\'ementaire de $F_f$. Ces deux ensembles jouent le
  r\^ole crucial dans l'\'etude de la dynamique de $f$. 
L'ensemble $J_f$  est le plus petit ensemble infini, ferm\'e et
  {\it compl\`etement invariant}, i.e. $f^{\pm 1}(J_f)=J_f$. 
Par cons\'equent, soit il 
  est \'egal \`a $\mathbb{P}^1$, soit il est d'int\'erieur 
  vide. 
  De plus, aucun point de $J_f$ n'est 
  isol\'e \cite[4.2]{Beardon}.
\par 
  Un point $z\in
  \mathbb{P}^1$ est dit {\it p\'eriodique} s'il existe un entier
  positif $n$ tel que $f^n(z)=z$. Ce point est appel\'e {\it r\'epulsif}
  (resp. {\it attractif, indiff\'erent} et {\it indiff\'erent rationnel}) si la
  d\'eriv\'ee de $f^n$ en $z$ est de module strictement plus grand que
  1 (resp. de module strictement plus petit que 1, de module 1 et  
  une racine de l'unit\'e). D'apr\`es un th\'eor\`eme de Fatou, 
  l'ensemble des points p\'eriodiques non r\'epulsifs est fini
  \cite[p.210]{Beardon}. 
  L'ensemble des points p\'eriodiques r\'epulsifs est infini et
  dense dans $J_f$ \cite[p.70]{BriendDuval,Beardon}. Un point $z$ est
  appel\'e {\it pr\'ep\'eriodique} (resp. {\it pr\'ep\'eriodique
  r\'epulsif, ...}) s'il existe un $n\in\mathbb{N}$ tel
  que $f^n(z)$ soit p\'eriodique (resp. p\'eriodique r\'epulsif, ...).\par
Dans cet article, on s'int\'eresse \`a d\'eterminer 
les fonctions rationnelles poss\'edant
  le m\^eme ensemble de Julia. 
  On montre facilement que
  l'ensemble de Julia de $f^k$ est \'egal \`a $J_f$ pour tout $k\geq
  1$. Dans \cite{Fatou,Julia}, Fatou et Julia ont montr\'e 
  que si $f$ et $g$ sont permutables
  ($fg=gf$), alors $J_f=J_g$ ({\it voir} \'egalement \cite[4.2.9]{Beardon}). 
  Il est facile de montrer que dans ce
  cas $f$ et $g$ admettent une infinit\'e de
  points p\'eriodiques communs. Ces deux
  propri\'et\'es restent valables pour les endomorphismes holomorphes
  permutables  
  de $\mathbb{P}^k$ \cite{HuYang, DinhSibony}.
  Si $f$, $g$ sont permutables et poss\`edent
  deux suites d'it\'er\'es disjointes alors l'une des conditions
  suivantes est vraie \cite{Ritt, Eremenko}:
\begin{enumerate} 
\item $J_f=J_g=\mathbb{P}^1$;
\item $f$ et $g$ sont conjugu\'ees simultanement aux applications du
  type $\lambda z^d$, i.e. il existe un automorphisme holomorphe
  $\sigma$ de $\mathbb{P}^1$ tel que $\sigma
  f\sigma^{-1}=\lambda_1z^{d_1}$ et  $\sigma g\sigma^{-1}=\lambda_2z^{d_2}$
  o\`u $d$, $d_1$, $d_2$ sont des entiers et $\lambda$, $\lambda_1$ et
  $\lambda_2$ sont des racines de l'unit\'e.
\item $f$ et $g$ sont conjugu\'ees simultanement aux applications du
  type $\pm \mbox{T}_d$ o\`u $\mbox{T}_d$ 
  est le {\it polyn\^ome de Tchebychev} de
  degr\'e $d\geq 2$ d\'efini par $\mbox{T}_d(\cos x):=\cos (dx)$. 
\end{enumerate}    
L'ensemble de Julia de l'application $\lambda z^d$ pour $|\lambda|=1$
et $|d|\geq 2$ est le cercle unit\'e. L'ensemble de Julia 
des applications $\pm\mbox{T}_d$ pour
$d\geq 2$ est le segment $[-1,1]$. Tous les couples de polyn\^omes 
dont les ensembles de Julia se co\"{\i}ncident,
sont d\'etermin\'es ({\it voir} par exemple 
\cite{BakerEremenko}, \cite{SchmidtSteinmetz}). Pour les fonctions
rationnelles, ce probl\`eme est r\'esolu par Levin et Przytycki 
dans le cas o\`u aucune
composante p\'eriodique de l'ensemble
de Fatou n'est un domaine parabolique, un disque de
Siegel ou un annneau de Herman \cite{LevinPrzytycki}. 
\par
Notons ${\cal C}_f$ l'ensemble des points critiques de $f$. Notre
r\'esultat principal est le th\'eor\`eme suivant:
\begin{theoreme} Soient $f$, $g$ deux fonctions rationnelles
  d\'efinies sur $\mathbb{P}^1$ de degr\'es $\geq 2$ telles que $J_f=J_g$. 
  Supposons que le c\^one tangent de $J_f$ en un de ses points est une
  r\'eunion finie de demi-droites et que tout point de
  ${\cal C}_f\cap J_f$ est pr\'ep\'eriodique r\'epulsif. Alors l'une des
  conditions suivantes est vraie:
\begin{enumerate}
\item $J_f=J_g$ est un cercle ou un arc d'un cercle pour une certaine
  coordonn\'ee de $\mathbb{P}^1$;
\item il existe $m_1,\ldots,m_k\in \mathbb{N}$ et $m\in \mathbb{N}^+$ 
  v\'erifiant $f^{m_1}g\ldots f^{m_k}g=f^m$. 
\end{enumerate} 
\end{theoreme}
\begin{exemple} Soient $h$ une fonction rationnelle, $p$ 
  un entier positif et $\alpha\not=1$ une $p$-i\`eme racine de
  l'unit\'e. Posons $f(z):=h(z^p)$ et $g:=\alpha h(z^p)$. Il est clair que
  $J_f=J_g$ et que les suites d'it\'er\'es de $f$ et de $g$ sont
  disjointes. Mais on a $fg=f^2\not =gf$. On peut trouver facilement une telle
  fonction $f$ poss\'edant des points p\'eriodiques indiff\'erents rationnels.
\end{exemple}
\section{D\'emonstration du th\'eor\`eme principal}
Soit $f$ une fonction holomorphe non constante 
d\'efinie dans un voisinage $U$ de
$0\in\mathbb{C}$ telle que $f(0)=0$. Un sous-ensemble $J\subset U$ est
appel\'e {\it invariant} (resp. {\it 
  compl\`etement invariant}) par $f$ au voisinage de $0$ 
s'il existe un voisinage $V\subset U$ de
$0$ tel que $f(J)\cap V=J\cap V$ (resp. $f^{\pm 1}(J)\cap V=J\cap
V$). Il est clair que $J$ est  
compl\`etement invariant si et seulement si 
$f^{\pm 1}(J)\cap V\subset J$.\par
Un ensemble $J\subset \mathbb{C}$ est dit {\it analytiquement lamin\'e}
en $z$ si $J\cap V=\sigma^{-1}(K)$ o\`u $V\subset
\mathbb{C}$ est un voisinage de $z$, $K$ est un sous-ensemble de
$\mathbb{R}$ et $\sigma:V\rightarrow \mathbb{R}$ 
est une fonction r\'eelle analytique de rang
maximal en tout point. L'ensemble $J$ est {\it analytiquement lamin\'e} dans
un ouvert  $U$ s'il est en tout point de $U$. On montre facilement que
$J$ poss\`ede deux laminations analytiques diff\'erentes en $z\in J$
si et seulement si $z$ appartient \`a l'int\'erieur de $J$.
\begin{proposition} Soient $f$, $g$ deux fonctions holomorphes
  inversibles d\'efinies dans un voisinage $U$ de $0$ telles que
  $f(0)=g(0)=0$ et $|f'(0)|\not = 1$. Soit $J\subset U$ un ferm\'e
  compl\`etement invariant par $f$ et par $g$ au voisinage de $0$. Alors l'une
  des deux conditions suivantes est vraie:
\begin{enumerate}
\item il existe un voisinage $V\subset U$ de $0$ tel que $J$ soit
  analytiquement lamin\'e dans $V\setminus\{0\}$;
\item il existe $(m,n)\in\mathbb{Z}^2-(0,0)$ tel que $f^m=g^n$. 
\end{enumerate} 
\end{proposition}
Il existe une application holomorphe 
$\varphi$, dite {\it de Poincar\'e}, d\'efinie
dans un petit disque $D$ centr\'e en $0$, telle que
$\varphi(0)=0$, $\varphi'(0)=1$ et $\varphi^{-1}f\varphi(t)=\lambda
t$ pour tout $t\in D$, 
o\`u $\lambda:=f'(0)\not=0$ est de module diff\'erent de $1$. 
Sans perdre en g\'en\'eralit\'e, on suppose,
pour la suite de la preuve de la proposition 1, que $f$ est lin\'eaire:
$f(t)=\lambda t$. Posons $\lambda':=g'(0)$ et fixons un disque
$D\subset\subset U$
suffisamment petit de rayon $r_0>0$. On peut aussi supposer 
que $|\lambda|>1$ car le cas contraire sera
  trait\'e de m\^eme mani\`ere en rempla\c cant $f$ par $f^{-1}$.
\begin{lemme}
Pour tout $t\in J\cap D$, pour tous entiers $m$, $n$ v\'erifiant
$\lambda^m{\lambda'}^nt\in D$, on a $\lambda^m{\lambda'}^nt\in J$.
\end{lemme}  
\begin{preuve} 
Soit $r>0$ suffisamment petit tel que $g^n(D_r)\subset D$ o\`u
$D_r:=\{t\in\mathbb{C}:\ |t|<r\}$. Posons $t'=\lambda^{-M}t$ pour $M$
suffisamment grand de sorte que $|t'|<r$ et ${\lambda'}^nt'\in D$. On a
$t'\in J$ car $J$ est invariant par $f^{-1}$ au voisinage de $0$. Posons
$t_k:=f^k g^n f^{-k}(t')$. Comme $t'\in J$, $t_k$
appartient \`a  $J$ lorsqu'il appartient \`a $D$.
On \'ecrit $g^n$ sous la forme d'une s\'erie de Taylor convergente:
$$g^n(t)={\lambda'}^nt+a_2t^2+a_3t^3+\cdots$$
On a:
$$f^k g^n f^{-k}(t')={\lambda'}^nt'+ 
a_2\lambda^{-k}{t'}^2+a_3\lambda^{-2k}  {t'}^3+\cdots$$
Par cons\'equent, $\lim t_k={\lambda'}^nt'\in D$. 
D'autre part, 
$J\cap D$ est ferm\'e dans $D$. On conclut que 
${\lambda'}^nt'\in J$ ou encore $\lambda^{-M}{\lambda'}^nt\in J$. 
Finalement,  
$\lambda^m{\lambda'}^nt$ appartient \`a $J$ lorsqu'il appartient \`a $D$
car $J$ est invariant par $f^{\pm 1}$.
\end{preuve}
\begin{lemme} Supposons que $\lambda^m{\lambda'}^n\not=1$ pour
  tout $(m,n)\in\mathbb{Z}^2-(0,0)$. Alors la condition 1 de la
  proposition 1 est vraie.
\end{lemme}
\begin{preuve} Comme $\lambda^{m}{\lambda'}^{n}\not =1$ pour tout
  $(m,n)\in\mathbb{Z}^2-(0,0)$, le groupe multiplicatif 
  ferm\'e engendr\'e par $\lambda$ et par $\lambda'$ contient un
  sous-groupe d'un param\`etre r\'eel $\{\lambda_r:=\exp(ar):\
  r\in\mathbb{R}\}$ o\`u $a\not =0$ est un nombre complexe avec
  $\mbox{Re}(a)\geq 0$. D'apr\`es le lemme
  pr\'ec\'edent, $J\cap D$ est une r\'eunion de courbes r\'eelles
  analytiques du type
  $\{t\in D:\ t=\lambda_r t_0\}$. Alors la condition 1 de la
  proposition 1 est vraie.
\par
Si $a$ est pure imaginaire, $J\cap D$ est une r\'eunion de
  cercles centr\'es en $0$.
\par
Si $a$ est un nombre r\'eel, $J\cap D$ est une r\'eunion de
  rayons de $D$.
\par
Sinon, $J\cap D$ est une
r\'eunion de courbes spirales tendant vers $0$.
\end{preuve}
D'apr\`es le lemme pr\'ec\'edent, il suffit de consid\'erer le cas
o\`u $\lambda^m={\lambda'}^n$ pour certain $(m,n)\not=(0,0)$. On a 
$n\not=0$ car $|\lambda|>1$. Posons $h:=f^mg^{-n}$. Supposons que la
condition 2 de la proposition 1 est fausse. On peut \'ecrire:
$$h(t)=t+\alpha t^{p+1}+\mbox{O}(t^{p+2})$$
o\`u $\alpha\not= 0$ et $p\geq 1$. Quitte \`a un changement de
coordonn\'ee du type $t\mapsto \sqrt[p]{-\alpha}t$, on peut
supposer que $\alpha=-1$. \par
Posons 
$$S_k:=\left\{t\in\mathbb{C}^*:\frac{(2k-1)\pi}{p}<\arg
    (t)<\frac{(2k+1)\pi}{p}\right\}$$ 
$$S_k(\epsilon):=\left\{t\in\mathbb{C}^*:\frac{(2k-1)\pi}{p}+\epsilon<\arg
    (t)<\frac{(2k+1)\pi}{p}-\epsilon\right\}$$
$$S(\epsilon):=\{t\in\mathbb{C}^*:-\pi+p\epsilon<\arg
    (t)<\pi-p\epsilon\}$$
et $$S:=\{t\in\mathbb{C}^*:\arg(t)\not =\pi\}$$
pour tout $0\leq k\leq p-1$ et tout $\epsilon>0$ suffisamment
    petit.\par
Soient $\Phi:\ \mathbb{C}^*\rightarrow\mathbb{C}^*$, $\Phi(t):=t^{-p}$ et
    $\Psi_k:\ S\rightarrow S_k$ les branches inverses de $\Phi$. Posons
    $r_1:=r_0^{-p}$, $S^r:=S\cap\{t:\ |t|>r\}$ pour $r>0$,
    $\Lambda:=\Phi f\Psi_k$, $h_k:=\Phi h \Psi_k$,
    $J_k:=\Phi(J\cap S_k\cap D)$. Alors 
    $\Phi(D\cap S_k)=S^{r_1}$, $J_k$ est un ferm\'e de 
    $S^{r_1}$, $\Lambda(w)=\lambda^{-p}w$ est ind\'ependant de $k$, 
    $\Lambda^{-n}(J_k\cap
    S(\epsilon))\subset J_k$,
    $\Lambda^n(J_k\cap S(\epsilon))\cap S^{r_1}\subset J_k$ lorsque
    $|\arg(\lambda^n)|<\epsilon$. Pour tout
    $\epsilon>0$ suffisamment petit, il existe $R_\epsilon>r_1$ tel que
    $h_k$ soit d\'efinie sur $S^R\cap S(\epsilon)$ et  
    $h_k(w)=w+p+\mbox{O}(|w|^{-1/p})$. On a \'egalement $h_k(J_k\cap
    S^{R_\epsilon}\cap S(\epsilon))\subset J_k$.
\begin{lemme} Fixons un $\epsilon$ suffisamment petit et un
    $R:=R_\epsilon$ suffisamment grand. 
    Si $w_0\in J_k$ est un point v\'erifiant $w_0+c\in
    S^{R_\epsilon}$ pour tout $c>0$, alors $w_0+c\in J_k$ pour tout $c>0$.  
\end{lemme}
\begin{preuve} Soit $\Theta\subset S^R$ 
un c\^one ferm\'e de sommet $w_0$, d'angle
    $0<2\delta<\pi$, dirig\'e par la
    demi-droite $\{w_0+c \mbox{ avec } c>0\}$. 
    Soit $C>0$ une constante v\'erifiant:
    $$\left|h_k(w)-w-p\right|<C|w|^{-1/p} \mbox{ sur }
    S^R\cap S(\epsilon).$$
On choisit une suite $\{n_j\}_{j=1}^\infty$ d'entiers positifs
    v\'erifiant:
\begin{enumerate}
\item[a.] $\lim n_j=+\infty$, $\lim \arg(\lambda^{n_j})=0$;
\item[b.] $\Lambda^{-n_j}(\Theta)\subset S^R\cap S(\epsilon)$ pour tout $j$;
\item[c.] $v_j:=|\lambda^{n_j}|>4CR^{-1/p}p^{-1}\delta^{-1}$ et 
      $0\leq 
\theta_j:=\arg(\lambda^{n_j})<\min(\epsilon,\delta/4p)$ pour tout $j$.
\end{enumerate}
Montrons d'abord que pour tout $j\geq 1$, on a $\Lambda^{n_j} h_k
\Lambda^{-n_j} (\Theta)\subset\Theta$. Soient $s\in \Theta$, 
$s_1:=\Lambda^{-n_j}(s)=\lambda^{n_jp}s$ et $s_2:=\Lambda^{n_j} h_k
\Lambda^{-n_j}(s)$. D'apr\`es b., $s_1\in S^R\cap S(\epsilon)$. Ceci
implique que $s_2$ est bien d\'efini. On a:
\begin{eqnarray*}
s_2 & = & \lambda^{-n_jp}h_k(\lambda^{n_jp}s)
\end{eqnarray*}
et
\begin{eqnarray*}
s_2-s & = &  \lambda^{-n_jp}(p+\gamma_s)
\end{eqnarray*}
avec $|\gamma_s|<C|\lambda|^{-n_j}|s|^{-1/p}<Cv_j^{-1}R^{-1/p}$ car
$s\in S^R$.\\
D'apr\`es c., on a:
\begin{eqnarray*}
|\arg(s_2-s)|
& \leq &
|\arg(\lambda^{-n_jp})|+|\arg(p+iCv_j^{-1}R^{-1/p})| \\
& \leq &
\delta/4+|\arg(p+ip\delta/4)|\\
& \leq & \delta/2
\end{eqnarray*}
D'o\`u $s_2\in\Theta$.
\par
On pose $m_j$ la partie enti\`ere de
$cp^{-1}v_j^p$
et $w_j:=[\Lambda^{n_j} h_k
\Lambda^{-n_j}]^{m_j}(w_0)$ pour tout $j\geq 1$. L'estimation
pr\'ec\'edente de $\gamma_s$ nous donne:
\begin{eqnarray*}
|w_j-w_0-c| & \leq
&|m_j\lambda^{-n_jp}p-c|+m_j|\lambda|^{-n_j p}C|v^{-1}R^{-1/p}\\
& \leq &
|m_jv_j^{-p}p\exp(-ip\theta_j)-c|+m_jCv_j^{-(p+1)}R^{-1/p}\\
& \leq &
|m_jv_j^{-p}p\exp(-ip\theta_j)-c\exp(-ip\theta_j)|+\\
& & \hspace{2cm}
+|c\exp(-ip\theta_j)-c|+cCp^{-1}R^{-1/p}v_j^{-1}\\
& \leq &
v_j^{-p}p|m_j-cp^{-1}v_j^p|+ 
c|\exp(-ip\theta_j)-1|+cCp^{-1}R^{-1/p}v_j^{-1}\\
& \leq &
v_j^{-p}p+ 
c|\exp(-ip\theta_j)-1|+cCp^{-1}R^{-1/p}v_j^{-1}
\end{eqnarray*}
D'apr\`es a., ce dernier terme tend vers $0$ quand $j$ tend vers
l'infini. Autrement dit, $w_j$ tend vers $w_0+c$. 
Comme $w_0\in J_k$, on a $w_j\in J_k$. Comme $J_k$ est
ferm\'e dans $S^{r_1}$, on a $w_0+c\in J_k$.
\end{preuve}
On choisit un $R'>>R$ de sorte que $w+c\in S^R$ pour 
tout $c\geq 0$ et tout 
$w\in S':=S(\epsilon)\cap S^{R'}$. Le lemme pr\'ec\'edent implique que
$J$ est analytiquement lamin\'e dans les $\Psi_k(S')$. En appliquant
le lemme pr\'ec\'edent aux fonctions $f$, $h^{-1}$ et \`a la
coordonn\'ee $t':=\sqrt[p]{-1}t$, on d\'eduit que $J$ est
analytiquement lamin\'e dans les $\sqrt[p]{-1}\Psi_k(S')$. Lorsque
$\epsilon$ est suffisamment petit, les ouverts $\Psi_k(S')$ et
$\sqrt[p]{-1}\Psi_k(S')$ recouvent un voisinage de $0$ sauf
le point $0$. Ceci termine la preuve de la
proposition 1. 
\par\hfill $\square$
\begin{proposition} Soient 
$$f(z)=z+\alpha z^{p+1}+{\rm O}(z^{p+2})$$
et
$$g(z)=z+\beta z^{q+1}+{\rm O}(z^{q+2})$$ 
deux fonctions holomorphes
  d\'efinies dans un voisinage $U$ 
  de $0$ o\`u $\alpha\not=0$, $\beta\not=0$ et
  $1\leq p <q$. Soit $J\subset U$ un ferm\'e 
  compl\`etement invariant par $f$ et par
  $g$ au voisinage de $0$. Alors il existe un voisinage $V\subset U$ de $0$ 
tel que $J$ soit analytiquement lamin\'e dans $V\setminus\{0\}$.
\end{proposition}
Soient $0\leq k\leq p-1$ un entier et $a>0$ un nombre r\'eel. On
d\'efinit {\it les p\'etales} $\Pi_k(a)$ et $\Pi'_k(a)$ par 
$$\Pi_k(a):=\{ r\exp(i\theta):\ 0<r^p<a(1+\cos(p\theta));
|2k\pi/p-\theta|<\pi/p\}$$
et
$$\Pi'_k(a):=\{ r\exp(i\theta):\ 0<r^p<a(1-\cos(p\theta));
|(2k+1)\pi/p-\theta|<\pi/p\}.$$ 
Alors pour tout $a>0$ les p\'etales $\Pi_k(a)$ et $\Pi'_k(a)$
recouvrent un voisinage de $0$ sauf le point
$0$. Les applications $\Phi$ et $\Psi_k$ sont d\'efinies dans la
preuve de la proposition 1. 
L'image $\Pi(a)$ de $\Pi_k(a)$ par $\Phi$
est ind\'ependante de $k$:
$$\Pi(a):=\{x+iy:\ y^2>1/a^2-2x/a\}.$$
D'apr\`es \cite[6.5.7]{Beardon}, il existe une coordonn\'ee
locale $t$ dans laquelle on a: 
$$f(t)=t-t^{p+1}+\mbox{O}(t^{2p+1})$$
et 
$$g(t)=t+\gamma t^{q+1}+\mbox{O}(t^{q+2}).$$
\begin{lemme}[{\cite[pp.116-122]{Beardon}}] 
  Soient $f(t)=t-t^{p+1}+\mbox{\rm O}(t^{2p+1})$ une fonction
  holomorphe d\'efinie au voisinage de $0$. Alors pour tout $k$
  et tout $a>0$
  suffisamment petit:
\begin{enumerate}
\item l'image de $\Pi_k(a)$ par $f$ est incluse dans $\Pi_k(a)$;
\item $f^n(w)\rightarrow 0$ localement uniform\'ement sur $\Pi_k(a)$;
\item $\arg f^n(w)\rightarrow 2k\pi/p$ localement uniform\'ement sur
  $\Pi_k(a)$; 
\item 
  il existe une application holomorphe injective
  $u:\Pi(a)\rightarrow \Pi(a)$ telle que $u(w)-w$ soit born\'ee,
  $\lim_{|w|\rightarrow\infty} 
  u'(w)=1$ et telle que
  $uvu^{-1}=T$ sur $u(\Pi(a))$ 
  o\`u $v:=\Phi f \Psi_k$ et $T(w):=w+p$.
\end{enumerate} 
\end{lemme}
Fixons un $a$ suffisamment petit. Posons $\Pi':=u(\Pi(a))$,
$w^{1/p}:=\Psi_0(w)$, $l:=\Phi g\Psi_k$, $\tau:=ulu^{-1}$,  
$I:=\Phi(J\cap\Pi_k(a))$ et $I':=u(I)$. D'apr\`es le lemme
pr\'ec\'edent, pour tout $a'$ suffisamment
petit, on a $\Pi(a')\subset \Pi'$. Comme $a$ est petit, $\tau$ et
$\tau^{-1}$ sont d\'efinies sur $\Pi'$. L'ensemble $I'$ est un ferm\'e 
{\it compl\`etement
invariant}  par $T$ et par $\tau$ dans le sens suivant:
$$T^{\pm 1}(I')\cap \Pi'\subset I'\mbox{ et }
\tau^{\pm 1}(I')\cap \Pi'\subset I'.$$ 
On peut \'ecrire $l$ sous la forme:
$$l(w)=w+\theta w^{-\delta/p}+\mbox{O}(w^{-(\delta+1)/p})$$
o\`u $\delta:=q-p\geq 1$ et $\theta\not=0$.
\\
Selon les propri\'et\'es de $u$ et sa d\'eriv\'ee ({\it voir} le lemme
  4), on a 
\begin{eqnarray*}
\tau(w) & = &
u\left(u^{-1}(w)+\theta[u^{-1}(w)]^{-\delta/p}+ 
\mbox{O}([u^{-1}(w)]^{-(\delta+1)/p})\right)\\
& = &
w+[\theta+\mbox{o}(1)][u^{-1}(w)]^{-\delta/p}+
\mbox{O}(w^{-(\delta+1)/p})\\
& = & w+\theta w^{-\delta/p}+ e(w)
\end{eqnarray*}
o\`u $e(w)=\mbox{o}(w^{-\delta/p})$.
\begin{lemme} Pour tout $w_0\in I'$ et tout  $r\in\mathbb{R}$, si
  $w_0+\theta r$ appartient \`a $\Pi'$, il appartient \`a $I'$.
\end{lemme}
\begin{preuve} Fixons un $r$, un $\epsilon>0$ et un
  $j\in\mathbb{N}$ v\'erifiant:
\begin{enumerate}
\item $B\subset \Pi'$ o\`u $B$ est la boule ferm\'ee
  de rayon $4|r|(|\theta|+\epsilon)$, de centre $x_0:=w_0+jp$;
\item $|e(w)|<\epsilon |w|^{-\delta/p}$ pour tout $w\in B$;
\item $|w|^{-\delta/p}<2(jp)^{-\delta/p}$ pour tout $w\in B$;
\item
  $|w^{-\delta/p}-(jp)^{-\delta/p}|<
8\delta[|r|(|\theta|+\epsilon)+|w_0|](jp)^{-1-\delta/m}$  
pour tout $w\in B$.
\end{enumerate} 
La condition 4 est r\'ealisable car la d\'eriv\'ee de la fonction
$w^{-\delta/p}$ est \'egale \`a $-\delta w^{-1-\delta/p}/p$ et car
$|w-jp|\leq
  4r(|\theta|+\epsilon)+|w_0|$.
\par
Par la suite, on traite le cas o\`u $r>0$; le cas contraire sera
trait\'e de m\^eme mani\`ere en utilisant la fonction $\tau^{-1}$ \`a
la place de $\tau$.
\par
Posons $M$ la partie enti\`ere de
  $(jp)^{\delta/p}r$ et $x_m:=\tau^{m}(x_0)$ pour tout
  $0\leq m\leq M$. On a pour tout $w\in B$:
\begin{eqnarray*}
|\tau(w)-w| &=& |\theta w^{-\delta/p}+e(w)|\\
& \leq & (|\theta|+\epsilon)|w|^{-\delta/p}\\
& \leq & 2(|\theta|+\epsilon)(jp)^{-\delta/p}
\end{eqnarray*}
ce qui implique par r\'ecurrence que les $x_m$ sont bien d\'efinis et
appartiennent \`a $B$ pour tout $1\leq m\leq M$. On a les estimations
suivantes:
\begin{eqnarray*}
|x_{m+1}-x_m-\theta(jp)^{-\delta/p}| & = & |\theta
x_m^{-\delta/p}-\theta(jp)^{-\delta/p}+ e(x_m)|\\
& \leq & |\theta
x_m^{-\delta/p}-\theta(jp)^{-\delta/p}| + |e(x_m)|\\
& \leq & |\theta||x_m^{-\delta/p}-(jp)^{-\delta/p}| + \epsilon
|x_m|^{-\delta/p}\\
& \leq & \{8\delta
|\theta|[r(|\theta|+\epsilon)+|w_0|](jp)^{-1}+2\epsilon\}(jp)^{-\delta/p}.
\end{eqnarray*}
D'o\`u: 
\begin{eqnarray*}
\lefteqn{|x_M-w_0-jp-M\theta(jp)^{-\delta/p}|=}\hspace{2.5cm}\\ 
& = &
|x_M-x_0-M\theta(jp)^{-\delta/p}|\\ 
& \leq & \sum_{m=0}^{M-1}|x_{m+1}-x_m-\theta(jp)^{-\delta/p}|\\
& \leq &
\{8\delta|\theta|[r(|\theta|+\epsilon)+|w_0|](jp)^{-1}
+2\epsilon\}(jp)^{-\delta/p}M\\
& \leq & \{8\delta|\theta|[r(|\theta|+\epsilon)+|w_0|](jp)^{-1}+2\epsilon\}r
\end{eqnarray*}
et 
\begin{eqnarray*}
|x_M-w_0-jp-\theta r| & \leq & |x_M-w_0-jp-M\theta(jp)^{-\delta/p}| +
|\theta r -M\theta(jp)^{-\delta/p}|\\
& \leq &  \{8\delta|\theta|[r(|\theta|+\epsilon)+|w_0|](jp)^{-1}+2\epsilon\}r +\\
& & \hspace{3.5cm} +|\theta|(jp)^{-\delta/p}|(jp)^{\delta/p}r-M|\\
& \leq & \{8\delta|\theta|[r(|\theta|+\epsilon)+|w_0|](jp)^{-1}+2\epsilon\}r+
|\theta|(jp)^{-\delta/p}.
\end{eqnarray*}
Comme $w_0\in I'$, on a $x_0\in I'$ et donc $x_M\in I'$ car $I'$ est
invariant par $T$ et par $\tau^{\pm 1}$. Pour $j$
tendant vers l'infini et $\epsilon$ tandant vers $0$, $x_M-jp$ tend vers
$w_0+\theta r$. Alors pour $j$ grand et $\epsilon$ petit  $x_M-jp\in I'$
car $I'$ est invariant par $T^{-1}$. La fermeture de $I'$
implique que $w_0+\theta r\in I'$.
\end{preuve}
D'apr\`es le lemme 5, $I'$ est
analytiquement lamin\'e en $w_0$. 
Par cons\'equent, $J$ est analytiquement lamin\'e dans les $\Pi_k(a')$ pour
$a'$ suffisamment petit. Maintenant, on remplace $f$ par $f^{-1}$ 
et la coordonn\'ee $t$ par $\sqrt[p]{-1}t$. Ces
remplacements ne changent pas les formes de $f$, $g$ et permutent les
deux familles de p\'etales $\{\Pi_k(a')\}$ et $\{\Pi'_k(a')\}$. On en
d\'eduit que $J$ est analytiquement lamin\'e dans les
$\Pi'_k(a)$. Ceci termine la preuve de la proposition 2 car 
les $\Pi_k(a')$ et $\Pi'_k(a')$ recouvrent un
voisinage de $0$ sauf le point $0$. 
\begin{corollaire} Soient $f$, $g$ deux fonctions holomorphes
  inversibles d\'efinies dans un voisinage $U$ de $0$ avec
  $f(0)=g(0)=0$ et $J\subset U$ un ferm\'e
  compl\`etement invariant par $f$ et par $g$ au voisinage de
  $0$. Supposons que le c\^one tangent de $J$ en $0$ n'est pas \'egal
  \`a $\mathbb{C}$. Alors l'une des conditions suivantes est vraie:
\begin{enumerate}
\item il existe un voisinage $V$ de $0$ tel que $J$ soit
  analytiquement lamin\'e dans $V\setminus\{0\}$;
\item il existe $(m,n)\in \mathbb{Z}^2-(0,0)$ tel que $f^m=g^n$.
\end{enumerate} 
\end{corollaire}
\begin{preuve}
Supposons que la condition 2 est fausse. D'apr\`es la proposition 1,
il suffit de traiter le cas o\`u $|f'(0)|=|g'(0)|=1$. Si $0$ n'est pas un
point d'accumulation de $J$, la condition 1 est vraie. Supposons que
ce n'est pas le cas. Comme le c\^one
tangent de $J$ en $0$ n'est ni vide ni \'egal \`a $\mathbb{C}$ et comme cet
ensemble est invariant par $f$ et par $g$, 
$f'(0)$ et $g'(0)$ sont des racines de l'unit\'e. Quitte \`a
remplacer $f$, $g$ par leurs it\'er\'es, on peut supposer que
$f'(0)=g'(0)=1$. On \'ecrit: 
$$f(z)=z+\alpha z^{p+1}+\mbox{O}(z^{p+2})$$ et 
$$g(z)=z+\beta z^{q+1}+\mbox{O}(z^{q+2})$$
o\`u $p\geq 1$ et $q\geq 1$ sont des entiers. Sans perdre en
g\'en\'eralit\'e, on peut supposer que $p\leq q$. En changeant la
coordonn\'ee, on peut supposer que $\alpha=-1$. D'apr\`es la
proposition 2, il nous reste \`a traiter le cas o\`u $p=q$. 
\begin{trivlist}
\item {\it Cas 1.} Supposons que $\beta\in\mathbb{Q}$. Quitte \`a
  remplacer $g$ par son it\'er\'e, on peut supposer que $\beta\in
  \mathbb{Z}$. Posons $h:=gf^{\beta}$. Comme la condition 2 est
  fausse, $h$ s'\'ecrit sous la
  forme: 
$$h(z)=z+\beta'z^{q'+1}+\mbox{O}(z^{q'+2})$$
avec $\beta'\not =0$ et $q'>p$.\\
En appliquant la proposition 2 aux fonctions $f$, $h$ et \`a
l'ensemble $J$, on constate que la condition 1 est vraie.
\item {\it Cas 2.} Supposons maintenant que $\beta\not \in
  \mathbb{Q}$. On utilise les notations de la preuve de la proposition
  2. On a:
$$f(t)=t-t^{p+1}+\mbox{O}(t^{2p+1})$$
$$g(t)=t+\beta t^{p+1}+\mbox{O}(t^{p+2})$$
et
$$\tau(w)=w-\beta p+\mbox{O}(w^{-1/p}).$$
\begin{lemme} Soit $w_0\in I'$.
\begin{enumerate}
\item Pour tous $m$, $n\in \mathbb{N}$ v\'erifiant $w_0+mp+np\beta\in
  I'$,  $w_0+mp+np\beta\in I'$.
\item Pour tout $r\in \mathbb{R}$ v\'erifiant
  $w_0+r\in \Pi'$, on a $w_0+r\in I'$.
\end{enumerate}
\end{lemme}
\begin{preuve} 1. Fixons un $\epsilon>0$. 
Choisissons un entier $s$ suffisamment
grand v\'erifiant:
\begin{enumerate}
\item $\Pi'':=\{w:\ {\rm Re\ } w\geq sp -|w_0|-|m|p-|n|p-\epsilon\}
\subset\Pi'$; 
\item $|\tau^{\pm 1}(w)- w \pm p\beta|<\epsilon/|n|$ pour tout
 $w\in \Pi''$.
\end{enumerate}
Posons $x_0:=w_0+sp+mp$ et 
$x_\nu:=\tau^{-{\rm sign}(n)}(x_{\nu-1})$ pour $\nu=1,\ldots, |n|$ o\`u ${\rm
  sign}(n)$ est le signe de $n$. 
Les $x_\nu$ sont bien d\'efinis et appartiennent \`a
$\Pi''$ car pour tout $1\leq \nu\leq |m|$:
$$
|s_{\nu}-s_{\nu-1}-{\rm sign}(n)p\beta|  =  
     |\tau^{-{\rm sign}(n)}(s_{\nu-1})-s_{\nu-1}-{\rm sign}(n) p\beta|
 \leq  \epsilon/|n|
$$
et 
\begin{eqnarray*}
|s_{\nu}-w_0-sp-mp-{\rm sign}(n)\nu p\beta| & = & |s_\nu-s_0-{\rm
  sign}(n)\nu p\beta|\\
& \leq & \sum_{k=1}^{\nu} |s_k-s_{k-1}-{\rm sign}(n)p\beta|\\ 
& \leq & \nu\epsilon/|n|.
\end{eqnarray*}
Les in\'egalit\'es pr\'ec\'edentes nous donnent:
$$ |s_{|n|}-w_0-sp-mp-np\beta| \leq \epsilon.$$  
On sait que $w_0\in I'$, $I'$ est compl\`etement invariant par $T$ et
par $\tau$. Donc $s_{|n|}\in I'$. Comme $I'$ est ferm\'e, en
consid\'erant $\epsilon\rightarrow 0$, on a $w_0+sp+mp+np\beta\in
I'$. L'invariance de $I'$ par $T^{-1}$ implique que $w_0+mp+np\beta\in
I'$. 
\par
2. Par hypoth\`ese, le c\^one tangent de $J$ en $0$ n'est pas \'egal
\`a $\mathbb{C}$. On d\'eduit de 1. que $\beta$ est r\'eel. D'apr\`es
1., $w_0+r\in I'$ car $\beta\not\in \mathbb{Q}$ et car $I'$ est ferm\'e.
\end{preuve}
Le lemme pr\'ec\'edent implique que $J$ est analytiquement lamin\'e
dans les $\Pi_k(a')$ pour $a'>0$ suffisamment petit. Comme dans la
proposition 2, on remplace $f$ par $f^{-1}$, $g$
par $g^{-1}$ et la coordonn\'ee $t$ par $\sqrt[p]{-1}t$. Ces
remplacements ne changent pas les formes de $f$, $g$ et permutent les
deux familles de p\'etales $\{\Pi_k(a')\}$ et $\{\Pi'_k(a')\}$. On en
d\'eduit que $I$ est analytiquement lamin\'e dans les $\Pi'_k(a')$. Ainsi,
ceci termine la preuve du corollaire.
\end{trivlist}
\end{preuve}
\begin{lemme} Soit $f$ une fonction rationnelle de degr\'e sup\'erieur
  ou \'egal \`a 2. 
\begin{enumerate}
\item[i.] Si le c\^one tangent de $J_f$ en un point p\'eriodique
  r\'epulsif $z_0$ est une r\'eunion finie de demi-droites, $J_f$ est inclus
  dans une courbe ou dans un arc r\'eel analytique ferm\'e.
\item[ii.]  Si dans un ouvert $U$ 
de $\mathbb{P}^1$, $J_f\cap U$ est un arc
  r\'eel de classe ${\cal C}^1$, alors $J_f$ est un
  cercle ou un arc d'un cercle pour une certaine coordonn\'ee.
\end{enumerate}
\end{lemme}
\begin{preuve} i. Quitte \`a remplacer $f$ par l'un de ses it\'er\'es, on peut
  supposer que $z_0$ est un point fixe. On peut \'egalement supposer
  que $z_0=0$. Soient $\varphi$ l'application de Poincar\'e de $f$ en
  $0$ et $\lambda:=f'(0)$. 
  On a $\varphi(0)=0$, $\varphi'(0)=1$ et $\varphi^{-1}
  f\varphi(t)=\lambda t$ au voisinage de $0$. Cette application
  $\varphi$ se prolonge en une application de $\mathbb{C}$ dans
  $\mathbb{P}^1$ par $\varphi(\lambda^n z):=f^n(\varphi(z))$ car
  $|\lambda|>1$. Posons
  $I:=\varphi^{-1}(J_f)$. Alors le c\^one tangent de $I$ en $0$ est
  une r\'eunion finie de demi-droites. Comme $J_f$ est invariant par
  $f^{-1}$, $I$ est invariant par $t\longmapsto\lambda^{-1} t$. Par
  cons\'equent, $\lambda$ est une racine d'un nombre r\'eel et $I$ est
  inclus dans une r\'eunion finie de droites. Quitte \`a remplacer
  $z_0$ par un autre point p\'eriodique pr\`es de $z_0$, on peut
  supposer que le c\^one tangent de $J_f$ en $z_0$ est une
  demi-droite ou une droite. Par cons\'equent, $I$ est inclus dans
  une droite passant par $0$ et $\lambda$ est r\'eel. 
Comme $J_f$ est r\'epulsif, il existe
  $z\not = 0$
  proche de $0$ et $n\in\mathbb{N}$ tel que $f^n(z)=0$. Il existe
  $x\in\mathbb{C}$ v\'erifiant $\varphi(x)=z$ car l'image de $\varphi$
  contient un voisinage de $0$. Posons $y:=\lambda^n x$. On a:
  $$\varphi(y)=\varphi(\lambda^n x)=f^n(\varphi(x))=f^n(z)=0.$$
Par cons\'equent, $J_f\subset \varphi([0,y])$, o\`u $[0,y]$ est
l'intervalle de sommets $0$ et $y$. L'ensemble $\varphi([0,y])$ est
une courbe ou un arc r\'eel analytique ferm\'e.
\par
ii. Comme l'ensemble des points p\'eriodiques r\'epulsifs de $f$
  est dense dans $J_f$, il existe un point $z_0\in J_f\cap U$ p\'eriodique
  r\'epulsif. D'apr\`es i., $I$ est inclus dans une droite
  r\'eelle $d$. Comme $J_f$ contient un arc r\'eel passant par $0$, $I$
  contient un voisinage de $0$ dans $d$. D'autre part, $I$ est
  invariant par $t\longmapsto\lambda t$ et $\lambda\in\mathbb{R}$,
  $|\lambda|>1$.  On en d\'eduit que $I$ est \'egal \`a $d$. 
Comme dans i., on montre que 
  $J_f=\varphi([0,y])$ qui est une courbe 
ou un arc r\'eel analytique ferm\'e. 
Cette courbe ou arc est lisse car $I=\varphi^{-1}(J_f)$ est lisse. Par
  cons\'equent, l'ensemble de Fatou $F_f$ contient une ou deux
  composantes dans lesquelle $f^n$ tend localement uniform\'ement vers
  un point fixe attractif ou indiff\'erent rationnel. D'apr\`es
  \cite{Fatou3}, $J_f$ est un cercle ou un arc d'un
  cercle pour une certaine coordonn\'ee.
\end{preuve}
\begin{proposition} Soient $f$ 
v\'erifiant l'hypoth\`ese du th\'eor\`eme 1 
et $z_0$ un point de $J_f$. 
  Supposons que $J_f$ n'est pas inclus dans une courbe ou dans un arc r\'eel
  analytique ferm\'e et que 
  le c\^one tangent de $J_f$ en $z_0$ est r\'eunion finie
  de demi-droites. Alors 
  $z_0$ est pr\'ep\'eriodique indiff\'erent rationnel. 
\end{proposition}
\begin{remarque} Dans \cite{BakerEremenko}, Baker et Eremenko ont
  d\'emontr\'e que $z_0$ est pr\'ep\'eriodique lorsque le c\^one
  tangent de $J_f$ en $z_0$ est une demi-droite et ${\cal C}_f\cap
  J_f=\emptyset$. Dans ce cas,
  l'hypoth\`ese ''$J_f$ n'est pas 
  inclus dans une courbe ou un arc r\'eel
  analytique'' n'est pas n\'ecessaire. Le r\'esultat de Baker et
  Eremenko
  reste valable
  sous l'hypoth\`ese de la proposition pr\'ec\'edente. La preuve dans
  \cite{BakerEremenko}  
  dont l'id\'ee principale a \'et\'e utilis\'ee dans \cite{Fatou3, Brolin},
  est essentiellement valable pour cette proposition.  
\end{remarque}
\begin{preuve}
La classification des domaines de Fatou en cinq types par Sullivin
\cite{Beardon}, tout point $z\in J_f$ d'accumulation de la suite
$\{f^n(c)\}_{n\geq 0}$ est p\'eriodique indiff\'erent rationnel 
pour tout $c\in F_f$ et en particulier pour $c\in {\cal C}_f\cap F_f$.
Notons
$A$ l'ensemble des points p\'eriodiques indiff\'erents rationnels de
$f$ et $B$ (resp. $C$) l'ensemble des points (resp. des points
p\'eriodiques) du type $f^n(c)$ avec $c\in {\cal C}_f\cap J$. 
Ces ensembles sont  finis. Quitte \`a
remplacer $f$ par un $f^m$, on peut supposer que tous les points de $A\cup C$
sont fixes. Soit $M>>0$ tel que $|f'(a)|<M$ pour tout $a\in A\cup C$.
\par
Supposons que $z_0$ n'est pas pr\'ep\'eriodique. Alors les
$f^n(z_0)$ sont tous diff\'erents. Il existe une suite
$\{n_\nu\}$ tendant vers l'infini et un $a\in J_f$ tels que
$\lim_{\nu\rightarrow \infty} f^{n_\nu}(z_0)=a$. Montrons qu'on peut
choisir $a\not\in A\cup B$. 
Supposons que ce n'est pas le cas. Alors on a $a\in A\cup B$. Quitte
\`a remplacer $n_\nu$ par $n_\nu+n$ pour $n>>0$, on peut supposer que
$a\in A\cup C$.  
On choisit, un disque $D(r)$ de centre $a$ et de rayon
$r>0$ suffisamment petit. On a $|z-a|<|f(z)-a|<2M|z-a|$
pour tout $z\in J_f\cap D(r)$. En effet, la deuxi\`eme in\'egalit\'e
est \'evidente, la premi\`ere pour $a$ r\'epulsif l'est aussi, si $a$
est indiff\'erent rationnel, ceci est \'enonc\'e dans
\cite[\S 31]{Fatou3}. Par cons\'equent, la courronne
$D(r)\setminus \overline{D(r/2M)}$ contient une infinit\'e de
$f^{n}(z_0)$. C'est contradiction car pour $r$ petit cette courronne ne
rencontre pas $A\cup C$.
\par
Notons $D$ un disque de centre $a$, de rayon $3\rho$ qui ne contient
aucun point de $\{f^n(c):\ c\in {\cal C}_f\}$ et
$z_\nu:=f^{n_\nu}(z_0)$. Soit $z=g_\nu(w)$ la
branche inverse de $w=f^{n_\nu}(z)$ telle que $g(z_\nu)=z_0$. Ces sont
des fonctions analytiques d\'efinies sur $D$.
D'apr\`es \cite[theorem 6.2]{Brolin}, 
la famille des $g_\nu$ est normale et pour tout
compact $\Delta\subset F_f\cap D$ ne contenant pas de points
p\'eriodiques, $g_\nu(\Delta)\rightarrow J_f$. Comme $J_f$ est
d'int\'erieur vide, on d\'eduit que toute fonction limite de la suite
$\{g_\nu\}$ est \'egale \`a la fonction constante $z_0$. On a aussi
$\lambda_\nu:= g'_\nu(z_\nu)\rightarrow 0$ quand $\nu\rightarrow
\infty$. Posons
$$\varphi_\nu(t):=[g_\nu(z_\nu+\rho t)-z_0]/(\rho\lambda_\nu)$$
qui appartiennent \`a la classe $S$ des fonctions holomorphes
univalentes d\'efinies sur $|t|<1$ normalis\'ees par $\varphi_\nu=0$,
$\varphi'_\nu(0)=1$. Remplacer $\{n_\nu\}$ par une sous-suite, on
peut supposer que $\arg\lambda_\nu\rightarrow \mu$ et
$\varphi_\nu\rightarrow \varphi\in S$ localement uniform\'ement dans
$|t|<1$. On obtient:
$$g_\nu(z_\nu+\rho
t)-z_0=\rho\lambda_\nu(\varphi(t)+\epsilon_\nu(t))$$
o\`u $\epsilon_\nu(t)\rightarrow 0$ localement uniform\'ement sur
$|t|<1$. 
\par
Pour tout $z\in J_f$ tel que $0<|z-z_0|<\rho/2$, posons $z_\nu+\rho
t_\nu=z$. Alors $t_\nu\rightarrow t:=(z-z_0)/\rho$. L'\'egalit\'e
pr\'ec\'edente implique:
$$g_\nu(z)-z_0=\rho\lambda_\nu(\varphi(t)+\delta_\nu)$$
o\`u $\varphi(t)\not=0$, $\delta_\nu\rightarrow 0$. Il existe une
direction tangente $\theta$ de $J_f$ en $z_0$ telle que
$\theta=\mu+\arg\varphi(t)$. Par cons\'equent, $z$ appartient \`a
la courbe analytique d\'efinie par
$\varphi((z-z_0)/\rho)=\tau\exp(i(\theta-\mu))$ o\`u $\tau\in
\mathbb{R}^+$. L'ensemble des $\theta$ est fini car le c\^one
tangent de $J_f$ en $z_0$ est
une r\'eunion finie de demi-droites. Rappelons que les points
p\'eriodiques r\'epulsifs de $f$ sont denses dans $J_f$.   
D'apr\`es le lemme 7, $J_f$ est inclus dans une courbe ou dans un arc
r\'eel analytique ferm\'e. C'est contradiction.
\par
Alors $z_0$ est un point pr\'ep\'eriodique. Il existe $n<m$ tels que
$f^n(z_0)=f^m(z_0)$. Le point $a:=f^n(z_0)$ est un point p\'eriodique de
$f$. D'apr\`es le lemme 7, 
$z_0$ est indiff\'erent car $J_f$ n'est pas inclus dans une 
courbe ou dans un arc r\'eel analytique. C'est un point indiff\'erent
rationnel car le c\^one tangent de $J_f$ en $0$ est une r\'eunion
finie de demi-droites. 
\end{preuve}
{\it Preuve du th\'eor\`eme 1.} Comme le c\^one tangent de $J_f$ en un
de ces points est une r\'eunion de demi-droites, $J_f$ est d'int\'erieur vide. 
Supposons que la condition 1 de ce th\'eor\`eme est fausse. D'apr\`es
le lemme 7, $J_f$ n'est pas une courbe ou un arc r\'eel analytique.
Lorsque $J_f$ n'est pas inclus
dans une courbe ou dans un arc r\'eel analytique ferm\'e, 
notons $E$ l'ensemble des points $z_0$ o\`u le c\^one tangent de $J_f$
est une r\'eunion finie de demi-droites. 
Lorsque $J_f$ est inclus dans une courbe ou un arc r\'eel
analytique, notons $E$ l'ensemble des points $z_0$ o\`u le c\^one
tangent de $J_f$ est une demi-droite. 
 D'apr\`es la
proposition 3 (pour le premier cas) et
la remarque 1 (pour le second cas), tout point de $E$ est
pr\'ep\'eriodique. Il est clair que $f(E)\subset E$ et $g(E)\subset
E$. Par hypoth\`ese, $E$ est non vide. Notons $A$
l'ensemble des points $z\in E$ p\'eriodiques pour $f$. 
\begin{lemme} L'ensemble $A$ est fini. 
\end{lemme}
\begin{preuve} Si $J_f$ n'est pas inclus dans une courbe ou dans un
  arc r\'eel analytique ferm\'e, 
d'apr\`es la proposition pr\'ec\'edente $A$
  est fini.\par
Supposons maintenant que $J_f$ est inclus dans une courbe ou dans un arc
  r\'eel analytique ferm\'e. D'apr\`es le lemme 7, $J_f$ est un ferm\'e
  totalement disconnexe de $L:=\varphi(\mathbb{R})$ qui est une courbe
  ou un arc r\'eel analytique ferm\'e
  ({\it voir} le lemme 7). 
Notons $]a_s,b_s[$ les composantes connexes de $L\setminus J_f$. Alors
  $E$ est l'ensemble des $a_s$, $b_s$ et $F_f$ est connexe. On a
  $f(L\setminus J_f)\subset L\setminus J_f$.\par
Si $F_f$ est le domaine d'attraction d'un point fixe attractif, il est
  clair que $f^n(]a_s,b_s[)$ passe par le point attractif pour
  tout $s$ et pour tout $n$ suffisamment grand. Il y a un nombre fini de
  $]a_s,b_s[$ qui passent par ce point attractif. Par cons\'equent,
  $A$ est fini.
\par
Sinon $f$ poss\`ede un point fixe indiff\'erent rationnel $a$. 
Dans ce cas, $|f(z-a)|>|z-a|$ pour $z\in
  J_f$ proche de $a$ et $f^n(K)$ tendent vers $a$ pour tout
  compact $K$ de $]a_s,b_s[$. On en d\'eduit que $a$ est le sommet 
de $f^n(]a_s,b_s[)$ pour tout $s$ et pour tout $n$ suffisamment
  grand. Il y a un nombre fini de $]a_s,b_s[$ qui poss\`edent $a_0$
  comme un sommet. Par cons\'equent, $A$ est fini.
\end{preuve}
Pour tout $a\in
A$, il existe $n_a\in\mathbb{N}$ tel que $f^{n_a}(a)\in A$. Comme $A$
est fini, il existe $z_0\in A$ et $n_1,\ldots, n_s$ tels que
$f^{n_1}gf^{n_2}g\ldots f^{n_s}g(z_0)=z_0$. Posons $h:=
f^{n_1}gf^{n_2}g\ldots f^{n_s}g$. Quitte \`a changer de coordonn\'ee,
on peut supposer que $z_0=0$. Supposons que la condition 2 du
th\'eor\`eme 1 est fausse. D'apr\`es le corollaire 1 (appliqu\'ee \`a
$f$, $h$ et $J_f$), l'ensemble $J_f$ est analytiquement lamin\'e dans
$V\setminus\{0\}$ o\`u $V$ est un voisinage de $0$. Selon les preuves
des lemmes 3, 5 et 6, pour tout $z\in
V\setminus\{0\}$, $J_f$ contient un arc r\'eel analytique joignant $z$
et $0$. Soient $z\in V\setminus\{0\}$ et $n>0$ tels que
$f^n(z)=0$. Alors $f^n(J_f\cap V_z)\subset J_f$ pour tout voisinage
$V_z$ de $z$. De plus, lorsque $V_z$ est petit, $J_f\cap V_z$ est
analytiquement lamin\'e. Les descriptions ci-dessus de $J_f$ au
voisinage de $0$ montrent que $J_f$ est une r\'eunion finie d'arcs
r\'eels analytiques. D'apr\`es le lemme 7, $J_f$ est un cercle ou un
arc d'un cercle. C'est une contradiction car la condition 1 du
th\'eor\`eme 1 est suppos\'ee fausse. 
\end{document}